\documentclass[pdflatex,sn-mathphys-num]{sn-jnl}

\usepackage{graphicx}%
\usepackage{multirow}%
\usepackage{amsmath,amssymb,amsfonts}%
\usepackage{amsthm}%
\usepackage{mathrsfs}%
\usepackage[title]{appendix}%
\usepackage{xcolor}%
\usepackage{textcomp}%
\usepackage{manyfoot}%
\usepackage{booktabs}%
\usepackage{algorithm}%
\usepackage{algorithmicx}%
\usepackage{algpseudocode}%
\usepackage{listings}%
\usepackage{tikz}%
\usepackage[all]{xy}%

\theoremstyle{thmstyleone}%
\newtheorem{theorem}{Theorem}
\newtheorem{proposition}[theorem]{Proposition}%
\newtheorem{lemma}[theorem]{Lemma}%
\theoremstyle{thmstyletwo}%

\theoremstyle{thmstylethree}%
\newtheorem{definition}{Definition}%

\begin{document}

\title[Cross-connections of the normed algebra of finite rank bounded operators on a Hilbert space]{Cross-connections of the normed algebra of finite rank bounded operators on a Hilbert space}


\author*[1]{\fnm{A.} \sur{Anju}}\email{anjuantony842@gmail.com}

\author[2]{\fnm{P. G.} \sur{Romeo}}\email{romeo\_parackal@yahoo.com}

\affil*[1]{\orgdiv{Department of Mathematics}, \orgname{Cochin University of Science and Technology}, \city{Ernakulam}, \postcode{682022}, \state{Kerala}, \country{India}}

\affil[2]{\orgdiv{Department of Mathematics}, \orgname{Cochin University of Science and Technology}, \city{Ernakulam}, \postcode{682022}, \state{Kerala}, \country{India}}


\abstract{This article examines the cross-connection between the normal normal category $\mathscr F(H)$ and the normal dual $\mathscr F(H^{'})$ of the finite-dimensional subspaces of a Hilbert space $H$ and it's dual space $H^{'}$. Further, we describe  the cross-connection semigroup, which is a normed algebra isomorphic to the normed algebra of finite rank bounded operators on $H$, then we characterize compact operators and their spectrum by the normal cones in the normal category of proper subspaces of $H$.}




\maketitle

\section{Introduction}\label{sec1}

Hall commenced the study of the structure of regular semigroups using partially ordered sets of principal left and right ideals of the semigroup \cite{h}. In 1974, Grillet abstractly characterized these partially ordered sets as regular partially ordered sets and established the relation between, which he referred to as a cross-connection of regular partially ordered sets \cite{g, r, i}. 
For each cross-connection, Grillet constructed a fundamental regular semigroup, known as the cross-connection semigroup. However a regular semigroup and its full regular subsemigroups have  isomorphic cross-connections, which reveals the limitation of the cross-connection of regular partially ordered sets. 
Subsequently, K. S. S. Nambooripad substituted regular partially ordered sets with normal categories and described the cross-connection of normal categories \cite{s, paa}. Nambooripad established that the cross-connection of normal categories produces a unique regular semigroup, called the cross-connection semigroup, and conversely, every regular semigroup is isomorphic to a cross-connection semigroup for an appropriate cross-connection \cite{s, pg}.

In this study, we consider the normal category $\mathscr F(H)$ of finite-dimensional subspaces of a Hilbert space $H$ with linear mappings as morphisms and investigate the relationship between the normal dual of $\mathscr F(H)$ and the dual space $H^{'}$ of $H$. 
It is shown that the regular semigroup of normal cones in $\mathscr F(H)$ is isomorphic to the semigroup of finite rank operators (linear maps) on $H$ \cite{a}. 
The bounded normal cones in $\mathscr F(H)$ are described and it is shown that these bounded normal cones in $\mathscr F(H)$ is a normed algebra isomorphic to the normed algebra $\mathcal S$ of all finite rank bounded operators on $H$ \cite{a}.
 Further, we prove that the normal dual of $\mathscr F(H)$ has a full subcategory, called the bounded normal dual of $\mathscr F(H)$, such that it is a normal category and is isomorphic to the normal category $\mathscr F(H^{'})$ of finite-dimensional subspaces of the dual space $H^{'}$ of $H$. Then the cross-connections of $\mathscr F(H^{'})$ to $\mathscr F(H)$ defined so that the cross-connection semigroup becomes a normed algebra isomorphic to the normed algebra $\mathcal S$ and we illustrate that for each invertible bounded operator $C$ on $H$, there exists a cross-connection $\Gamma_C$ from $\mathscr F(H^{'})$ to $\mathscr F(H)$ and a 
 cross-connection semigroup.  However, if $C$ is a unitary operator on $H$, then the cross-connection semigroup determined by $\Gamma_C$ is a normed algebra isomorphic to the normed algebra $\mathcal S$. 
 We also discuss the spectrum of compact operators on $H$ using normal cones in the normal category  $\mathscr P(H)$ of proper subspaces of $H$ with linear mappings as morphisms and 
 the Banach algebra $\mathcal K$ of all compact operators on an infinite-dimensional Hilbert space $H$,and is establish that the Banach algebra $\mathcal K$ is isomorphic to the closure of a subsemigroup of normal cones in $\mathscr P(H)$.,

\section{Normal categories and normal duals }\label{sec2}

Initially, we recall some of the terminology introduced by Nambooripad in the cross-connection theory in \cite{s}. Throughout this article, the order of the composition of functions, morphisms and functors is from left to right. For a category $\mathscr C$, we let $\textit{\textbf v} \mathscr C$ be the object class of $\mathscr C$, $\mathscr C$ itself be the morphism class of $\mathscr C$, and $1_{\mathscr C}$ be the identity functor on $\mathscr C$. 
For $a, b \in \textit{\textbf v} \mathscr C$, let $1_{a}$ be the identity morphism on $a$, and denote the collection of all morphisms from $a$ to $b$ by the hom-set $\mathscr C(a, b)$.

A category $\mathscr P$ is called a preorder if each hom-set of $\mathscr P$ contains at most one morphism. If $\mathscr P$ is a preorder, we can define a quasiorder $\subseteq$ on $\textit{\textbf v} \mathscr P$ as $p\subseteq p'$ if and only if $\mathscr P(p, p')\neq \varnothing$ for $p, p'\in \textit{\textbf v} \mathscr P$. If the relation $\subseteq$ on $\textit{\textbf v} \mathscr P$ is a partial order, then the preorder $\mathscr P$ is a strict preorder. A small category $\mathscr C$ is said to be a category with subobjects if there exists a subcategory $\mathscr P$ of $\mathscr C$, which is a strict preorder with $\textit{\textbf v} \mathscr P=\textit{\textbf v} \mathscr C$ such that every morphism in $\mathscr P$ is a monomorphism in $\mathscr C$ and if $f=hg$ for $f, g\in \mathscr P$ and $h\in \mathscr C$, then $h\in \mathscr P$. The morphisms in $\mathscr P$ are known as inclusions in $\mathscr C$, and if $q\in \mathscr C$ is a right inverse of an inclusion in $\mathscr C$, then $q$ is a retraction. For $f\in \mathscr C$, a factorization of the form $f=quj$ with $q$ a retraction, $u$ an isomorphism, and $j$  an inclusion is termed a normal factorization of $f$. In this case, the epimorphic component $f^{\circ}$ of $f$ is $qu$ \cite{s}.
\begin{definition}
Let $\mathscr C$ be a category with subobjects. A cone $\gamma$ in $\mathscr C$ with vertex $d\in \textit{\textbf v} \mathscr C$ is a map from $\textit{\textbf v} \mathscr C$ to $\mathscr C$ such that for $c, c' \in \textit{\textbf v} \mathscr C$, $\gamma(c)\in \mathscr C(c, d)$ and $\gamma(c')=j\gamma(c)$ if $j$ is an inclusion from $c'$ to $c$.
\end{definition}
The morphism $\gamma(c)$ is known as the component of $\gamma$ at $c$, and the vertex of $\gamma$ is denoted by $c_{\gamma}$. 
A cone $\gamma$ with vertex $d$ is a normal cone if there exists at least one $c\in \textit{\textbf v} \mathscr C$ such that $\gamma(c)$ is an isomorphism from $c$ to $d$. 
A category $\mathscr C$ with subobjects is referred to as a normal category if every morphism in $\mathscr C$ has a normal factorization, every inclusion in $\mathscr C$ splits and for each $d \in \textit{\textbf v} \mathscr C$, there is a normal cone $\gamma$ with vertex $d$ such that $\gamma(d)=1_d$. For a normal category $\mathscr C$, if we define the binary operation on the set $T\mathscr C$ of all normal cones in $\mathscr C$ by $(\gamma \cdot \delta)(c)= (\gamma\star(\delta(c_{\gamma}))^{\circ})(c)=\gamma(c)(\delta(c_{\gamma}))^{\circ}$ for $\gamma, \delta \in T\mathscr C$ and $c\in \textit{\textbf v} \mathscr C$, then $T\mathscr C$ is a regular semigroup \cite{s}.

If $\mathscr C$ is a normal category, then the normal dual $N^\ast \mathscr C$ of $\mathscr C$ is the full subcategory of the category $\mathscr C^{\ast}$ of all functors from $\mathscr C$ to $\mathbf{Set}$ such that 
$$\textit{\textbf v} N^\ast \mathscr C=\{H(\epsilon, -):\epsilon \in E(T\mathscr C)\},$$ 
where $E(T\mathscr C)$ is the set of all idempotents in $T\mathscr C$ and for $\epsilon \in E(T\mathscr C)$, $H(\epsilon, -):\mathscr C\rightarrow \mathbf{Set}$ is the $H$- functor \cite{s}. We have that $N^\ast \mathscr C$  is a normal category, and if $\mu: H(\epsilon, -)\rightarrow H(\epsilon', -)$ is a morphism in $N^*\mathscr C$, then there exists a unique morphism $\hat{\mu}: c_{\epsilon'}\rightarrow c_\epsilon$ in $\mathscr C$ such that for $c\in \textit{\textbf v} \mathscr C$ and $f\in \mathscr C(c_\epsilon, c)$,  the component of the natural transformation $\mu$ at $c$ is given by  
\begin{align*}
	\mu(c)&:H(\epsilon, c)\rightarrow H(\epsilon', c)\\
	&:\epsilon\star f^\circ \mapsto \epsilon'\star(\hat{\mu}f)^\circ.
\end{align*}

\subsection{Normal dual of $\mathscr F(H)$ and dual space of $H$}

Let $H$ be a Hilbert space over a field $\textbf K=\mathbb R$ or $\mathbb C$. Consider the category $\mathscr F(H)$ of finite-dimensional subspaces of $H$, whose objects are finite-dimensional subspaces of $H$, and for $M, N\in \textit{\textbf v} \mathscr F(H)$, the hom-set $\mathscr F(H)(M, N)$ is the set of all linear maps from $M$ to $N$. It is seen that $\mathscr F(H)$ with the usual subspace inclusions is a category with subobjects, and is a normal category. Further, the regular semigroup $T\mathscr F(H)$ of all normal cones in $\mathscr F(H)$ is isomorphic to the semigroup $\mathcal S^{'}$ of all finite rank operators on $H$ (see \cite[Theorem 5.1]{a}), also for each $\gamma\in T\mathscr F(H)$ with vertex $N$, there is a unique finite rank operator $T$ such that the range space $R(T)$ of $T$ is $N$ and $\gamma=\gamma^T:\textit{\textbf v} \mathscr F(H) \rightarrow \mathscr F(H)$ is given for $M\in \textit{\textbf v} \mathscr F(H)$ by
$$\gamma^T(M)=T|_M:M\rightarrow N.$$ 
Since $\gamma^{T_1} \cdot \gamma^{T_2}=\gamma^{T_1T_2}$ for $T_1, T_2\in \mathcal S^{'}$, $\gamma^{T}$ is an idempotent in $T\mathscr F(H)$ if and only if $T$ is an idempotent in $\mathcal S^{'}$.  The normal dual $N^\ast \mathscr F(H)$ of $\mathscr F(H)$ is the full subcategory of the category $\mathscr F(H)^\ast$ of all functors from $\mathscr F(H)$ to $\mathbf {Set}$ such that 
$$\textit{\textbf v} N^*\mathscr F(H)=\{H(\gamma^P, -): P\in E(\mathcal S^{'})\}.$$ 
For $P\in E(\mathcal S^{'})$ with $N=R(P)$, the $H$- functor $H(\gamma^P, -):\mathscr F(H) \rightarrow \mathbf {Set}$ is given for $M_1, M_2\in \textit{\textbf v} \mathscr F(H)$ and $T\in \mathscr F(H)(M_1, M_2)$ by
\begin{align*}
	H(\gamma^P, M_1)&=\{\gamma^P\star A^\circ : A\in \mathscr F(H)(N, M_1)\} \text{ and }\\
	H(\gamma^P, T)&:H(\gamma^P, M_1)\rightarrow H(\gamma^P, M_2)\\
	&:\gamma^P\star A^\circ \mapsto \gamma^P\star (AT)^\circ.
\end{align*}
However, for $A\in \mathscr F(H)(N, M_1)$, $\gamma^{PAP_{M_1}}=\gamma^P \cdot \gamma^{P_NAP_{M_1}}=\gamma^P \star ((P_NAP_{M_1})|_N)^\circ=\gamma^P \star A^\circ$, where $P_N$ and $P_{M_1}$ are projections (self-adjoint idempotents) of $H$ onto $N$ and $M_1$, respectively. Therefore, we have
\begin{align*}
	H(\gamma^P, M_1)&=\{\gamma^{PAP_{M_1}} : A\in \mathscr F(H)(N, M_1)\} \text{ and }\\
	H(\gamma^P, T)&:\gamma^{PAP_{M_1}} \mapsto \gamma^{PATP_{M_2}}.
\end{align*} 
Let $\mu:H(\gamma^{P_1}, -) \rightarrow H(\gamma^{P_2}, -)$ be a natural transformation in $N^*\mathscr F(H)$ for $P_1, P_2\in E(\mathcal S^{'})$, then a linear map $T:R(P_2)\rightarrow R(P_1)$ exists uniquely such that for $M\in \textit{\textbf v} \mathscr F(H)$ and $A\in \mathscr F(H)(R(P_1), M)$, 
\begin{align*}
	\mu(M)&:H(\gamma^{P_1}, M) \rightarrow H(\gamma^{P_2}, M)\\
	&:\gamma^{P_1AP_M}\mapsto \gamma^{P_2TAP_M}.
\end{align*}
We denote this natural transformation, determined by $T$, as $\mu_{T^{*}}$. Now define the bounded normal cones in $\mathscr F(H)$ as follows:
\begin{definition}\emph{\cite[Definition 5.1]{a}}
	A normal cone $\gamma$ in $\mathscr F(H)$ with vertex $N$ is called a bounded normal cone if there exists some $\alpha_\gamma\in \mathbb R$ such that $\lVert \gamma(M)\rVert\leq \alpha_\gamma$ for all $M\in \textit{\textbf v} \mathscr F(H)$.
\end{definition}
Let $B\mathscr F(H)$ denote the set of all bounded normal cones in $\mathscr F(H)$, due to 
Theorems 4.1 and 5.2 in \cite{a}, the set $B\mathscr F(H)$ of all bounded normal cones in $\mathscr F(H)$ is a normed algebra isomorphic to the normed algebra $\mathcal S$ of all finite rank bounded operators on $H$ under the operation given by, for $T, T_1, T_2 \in \mathcal S$ and $k\in \mathbf K$,
\begin{align*}
	\gamma^{T_1}+\gamma^{T_2}&=\gamma^{T_1+T_2},\\
	k\gamma^T&=\gamma^{kT},\\
	\lVert \gamma^T\rVert&=\lVert T\rVert \text{ and}\\
	\gamma^{T_1} \cdot \gamma^{T_2}&=\gamma^{T_1T_2}.
\end{align*}
Next, the bounded normal dual of $\mathscr F(H)$ is defined as  
\begin{definition}
	Let $\mathscr F(H)$ be the category of finite-dimensional subspaces of $H$. The bounded normal dual of $\mathscr F(H)$, denoted by $B^*\mathscr F(H)$, is defined as the full subcategory of the category $\mathscr F(H)^{\ast}$ of all functors from $\mathscr F(H)$ to $\mathbf{Set}$ with
	$$\textit{\textbf v} B^*\mathscr F(H)=\{H(\gamma^P, -): \gamma^P\in E(B\mathscr F(H))\}.$$ 
\end{definition}
It is apparent that $\gamma^P\in E(B\mathscr F(H))$ if and only if $P\in E(\mathcal S)$, and for  normal cones $\gamma$ and $\gamma'$ in a normal category $\mathscr C$, $H(\gamma, -)=H(\gamma', -)$ [$H(\gamma, -) \subseteq H(\gamma', -)$] if and only if there is a unique isomorphism [epimorphism] $h:c_{\gamma'}\rightarrow c_\gamma$ such that $\gamma=\gamma'\star h$ \cite[III, Proposition 7]{s}. Therefore, if $P\in E(\mathcal S)$ with $N=Z(P)^\perp$, then $P=P_NP$ implies $\gamma^P=\gamma^{P_N}\cdot\gamma^{P}=\gamma^{P_N}\star(P|_N)^\circ=\gamma^{P_N}\star P|_N$. Consequently, $H(\gamma^P, -)=H(\gamma^{P_N}, -)$, leads to the following lemma.
\begin{lemma}
	Let $\mathscr F(H)$ be the category of finite-dimensional subspaces of $H$. Then, the bounded normal dual $B^*\mathscr F(H)$ of $\mathscr F(H)$ is the full subcategory of $N^*\mathscr F(H)$ such that
$$\textit{\textbf v} B^*\mathscr F(H)=\{H(\gamma^{P_N}, -): N\in \textit{\textbf v} \mathscr F(H)\}.$$ 
\end{lemma}
Next, consider the normal category $\mathscr F(H^{'})$ of finite-dimensional subspaces of the dual space $H^{'}$ of the Hilbert space $H$. Then 
$$\textit{\textbf v} \mathscr F(H^{'})=\{M^{'}: M\in \textit{\textbf v} \mathscr F(H)\},$$
where $M^{'}=\{f_m\in H^{'}: m\in M\}$ and $f_m(x)=\langle x, m\rangle$ for $x\in H$, and for $M, N\in \textit{\textbf v} \mathscr F(H)$,
$$\mathscr F(H^{'})(M^{'}, N^{'})=\{(P_NT^{*}P_M)^{'}|_{M^{'}}: T\in \mathscr F(H)(M, N)\},$$
where $T^*$ is the adjoint of $T$ and  $(P_NT^*P_M)^{'}$ is the transpose of $P_NT^*P_M$, 
(see \cite{a}). Note that $(P_NT^{*}P_M)^{'}|_{M^{'}}$ maps $f_m\in M^{'}$ to $f_{T(m)}\in N^{'}$. Now define $F:\mathscr F(H^{'})\rightarrow B^*\mathscr F(H)$ by
$$F(N^{'})=H(\gamma^{P_N}, -)$$
for $N^{'}\in \textit{\textbf v} \mathscr F(H^{'})$ and for $(P_{N_2}TP_{N_1})^{'}|_{N_1^{'}}:N_1^{'}\rightarrow N_2^{'},$
$$F((P_{N_2}TP_{N_1})^{'}|_{N_1^{'}})=\mu_{T^{*}}:H(\gamma^{P_{N_1}}, -)\rightarrow H(\gamma^{P_{N_2}}, -).$$
Then, for $(P_{N_2}T_1P_{N_1})^{'}|_{N_1^{'}}:N_1^{'}\rightarrow N_2^{'}$ and $(P_{N_3}T_2P_{N_2})^{'}|_{N_2^{'}}:N_2^{'}\rightarrow N_3^{'},$
\begin{align*}
	F((P_{N_2}T_1P_{N_1})^{'}|_{N_1^{'}})(P_{N_3}T_2P_{N_2})^{'}|_{N_2^{'}}))&=F((P_{N_3}T_2T_1P_{N_1})^{'}|_{N_1^{'}})\\
	&=\mu_{(T_2T_1)^{*}}=\mu_{T_1^{*}}\mu_{T_2^{*}}\\
	&=F((P_{N_2}T_1P_{N_1})^{'}|_{N_1^{'}})F((P_{N_3}T_2P_{N_2})^{'}|_{N_2^{'}}).
\end{align*}
Additionally, if $M^{'}\subseteq N^{'}$, then $M\subseteq N$ implies that $(P_NP_M)^{'}|_{M^{'}}$ is the inclusion from $M^{'}$ to $N^{'}$, and $F((P_NP_M)^{'}|_{M^{'}})=\mu_{(P_M|_N)^*}$ is the inclusion from $H(\gamma^{P_M}, -)$ to $H(\gamma^{P_N}, -)$. That is, $F$ preserves identities and inclusions.

Moreover, if $H(\gamma^{P_N}, -)=H(\gamma^{P_M}, -)$, there exists a unique isomorphism $A$ from $M$ to $N$ such that $\gamma^{P_N}=\gamma^{P_M} \star A$. As a result, $\gamma^{P_N}(M)=(\gamma^{P_M} \star A)(M)$ implies $P_N|_M=A$. Since $P_N|_M$ is an isomorphism from $M$ to $N$, it follows that $\gamma^{P_N}=\gamma^{P_M} \cdot \gamma^{P_N}$. Thus, $M^\perp \subseteq N^\perp$ leads to $M=N$ and $M^{'}=N^{'}$. Therefore, $\textit{\textbf v}F:\textit{\textbf v}\mathscr F(H^{'})\rightarrow \textit{\textbf v}B^*\mathscr F(H)$ is a bijection, and for each natural transformation $\mu_{T^{*}}:H(\gamma^{P_N}, -)\rightarrow H(\gamma^{P_M}, -)$, there exists a unique linear map $T:M\rightarrow N$ such that $F((P_MTP_N)^{'}|_{N^{'}})=\mu_{T^{*}}$. Hence, $F$ is fully faithful. Also it is known that, if an inclusion preserving isomorphism exists between two normal categories, then they are isomorphic \cite{s}. Consequently, $B^*\mathscr F(H)$ is a normal category, and $F$ is a normal category isomorphism from $\mathscr F(H^{'})$ to $B^*\mathscr F(H)$. This indicates that the Hilbert space duality is equivalent to the bounded normal duality of $\mathscr F(H)$ and thus the following theorem.
\begin{theorem}
	Let $H$ be a Hilbert space and $\mathscr F(H)$ be the category of finite-dimensional subspaces of $H$. Then, the bounded normal dual $B^*\mathscr F(H)$ of $\mathscr F(H)$ is a normal category, and the normal category $B^*\mathscr F(H)$ is isomorphic to the normal category $\mathscr F(H^{'})$ of finite-dimensional subspaces of the dual space $H^{'}$ of $H$.
\end{theorem}

\section{Cross-connections}\label{sec3}

Now, we proceed to investigate the cross-connections of $\mathscr F(H^{'})$, the category of finite-dimensional subspaces of the dual Hilber space $H'$  of $H$ to $\mathscr F(H)$. To begin, we briefly recall cross-connection of normal categories as discussed in \cite{s}.

Let $\mathscr C$ be a normal category and $c\in \textit{\textbf v}\mathscr C$. The principal ideal generated by $c$, denoted by $\langle c \rangle$, is the full subcategory of $\mathscr C$ such that the objects of $\langle c \rangle$ are subobjects of $c$ in $\mathscr C$. For two normal categories $\mathscr C$ and $\mathscr D$, an inclusion preserving functor $F:\mathscr C\rightarrow \mathscr D$ is a local isomorphism if $F$ is fully faithful and $F|_{\langle c \rangle}$ is an isomorphism from $\langle c \rangle$ onto $\langle F(c) \rangle$ for every $c\in \textit{\textbf v}\mathscr C$. The $\mathcal M$-set $\mathcal MH(\gamma, -)$ of a normal cone $\gamma$ in $\mathscr C$ is defined by 
$$\mathcal MH(\gamma, -)=\{c\in \textit{\textbf v}\mathscr C: \gamma(c) \text{ is an isomorphism}\}.$$ For two normal categories $\mathscr C$ and $\mathscr D$, a triplet $(\mathscr D, \mathscr C; \Gamma)$ is a cross-connection from $\mathscr D$ to $\mathscr C$ if $\Gamma:\mathscr D\rightarrow N^*\mathscr C$ is a local isomorphism such that for every $c\in \textit{\textbf v}\mathscr C$, there exists some $d\in \textit{\textbf v}\mathscr D$ with $c\in \mathcal M\Gamma(d)$.
Also for each cross-connection $(\mathscr D, \mathscr C; \Gamma)$, there is a unique dual cross-connection $(\mathscr C, \mathscr D; \Delta)$ such that there exists a natural isomorphism $\chi_{\Gamma}$ from the bifunctor $\Gamma(-, -):\mathscr C\times \mathscr D\rightarrow \mathbf{Set}$ to the bifunctor $\Delta(-, -):\mathscr C\times \mathscr D\rightarrow \mathbf{Set}$, where $\chi_{\Gamma}$ is referred to as the duality associated with $\Gamma$. We can link the normal cones in $\mathscr C$ and $\mathscr D$ using the natural isomorphism $\chi_\Gamma$ in the following way: a cone $\gamma\in T\mathscr C$ is linked to a cone $\gamma^*\in T\mathscr D$ if there exists $(c, d)\in \textit{\textbf v}\mathscr C\times \textit{\textbf v}\mathscr D$ such that $\gamma \in \Gamma(c, d)$ and $\gamma^*= \chi_\Gamma(c, d)(\gamma).$ The set $\tilde{S}\Gamma$ of pairs of linked cones $(\gamma, \gamma^*)$ is a regular subsemigroup of $T\mathscr C\times (T\mathscr D)^{op}$, known as the cross-connection semigroup determined by $\Gamma$ \cite{s}.

Next, we proceed to find cross-connections of $\mathscr F(H^{'})$ to $\mathscr F(H)$. Let $C$ be an invertible bounded operator on $H$, and for $N\in \textit{\textbf v}\mathscr F(H)$, let $Q_{C^{-1}(N)}$ denote the idempotent operator on $H$ with $C^{-1}(N)$ as its range space and $C^{-1}(N^\perp)$ as its kernel. Since $C^{-1}(N)$ and $C^{-1}(N^\perp)$ are closed subspaces of $H$, $Q_{C^{-1}(N)}$ is a finite rank bounded operator on $H$. Define $\Gamma_C: \mathscr F(H^{'})\rightarrow N^{*}\mathscr F(H)$ by for $N^{'}\in \textit{\textbf v}\mathscr F(H^{'})$,
$$\Gamma_C(N^{'})=H(\gamma^{Q_{C^{-1}(N)}}, -) $$
and for $(P_{N_2}TP_{N_1})^{'}|_{N_1^{'}}:N_1^{'}\rightarrow N_2^{'}$,
$$\Gamma_C((P_{N_2}TP_{N_1})^{'}|_{N_1^{'}})=\mu_{((CTC^{-1})|_{C^{-1}(N_2)})^{*}}:H(\gamma^{Q_{C^{-1}(N_1)}}, -) \rightarrow H(\gamma^{Q_{C^{-1}(N_2)}}, -).$$
Then, for $(P_{N_2}T_1P_{N_1})^{'}|_{N_1^{'}}:N_1^{'}\rightarrow N_2^{'}$ and $(P_{N_3}T_2P_{N_2})^{'}|_{N_2^{'}}:N_2^{'}\rightarrow N_3^{'}$,
\begin{align*}
	\Gamma_C(((P_{N_2}T_1P_{N_1})^{'}|_{N_1^{'}})((P_{N_3}T_2P_{N_2})^{'}|_{N_2^{'}}))&=\Gamma_C((P_{N_3}T_2T_1P_{N_1})^{'}|_{N_1^{'}})\\
	&=\mu_{((CT_2T_1C^{-1})|_{C^{-1}(N_3)})^{*}}\\
	&=\mu_{((CT_1C^{-1})|_{C^{-1}(N_2)})^{*}}\,\mu_{((CT_2C^{-1})|_{C^{-1}(N_3)})^{*}}\\
	&=\Gamma_C((P_{N_2}T_1P_{N_1})^{'}|_{N_1^{'}})\Gamma_C((P_{N_3}T_2P_{N_2})^{'}|_{N_2^{'}}).
\end{align*}
Also, if $M^{'}\subseteq N^{'}$, then $(P_NP_M)^{'}|_{M^{'}}$ is the inclusion from $M^{'}$ to $N^{'}$ and $\Gamma_C((P_NP_M)^{'}|_{M^{'}})=\mu_{((C(P_M|_N)C^{-1})|_{C^{-1}(N)})^{*}}$ is the inclusion from $H(\gamma^{Q_{C^{-1}(M)}}, -)$ to $H(\gamma^{Q_{C^{-1}(N)}}, -)$. Hence, $\Gamma_C$ preserves identities and inclusions.  Further, for each natural transformation $\mu$ from $H(\gamma^{Q_{C^{-1}(N_1)}}, -)$ to $H(\gamma^{Q_{C^{-1}(N_2)}}, -)$, there exists a unique linear map $S$ from $C^{-1}(N_2)$ to $C^{-1}(N_1)$ such that $\mu=\mu_{S^{*}}$. Let $T=(C^{-1}SC)|_{N_2}$, then $T$ is a linear map from $N_2$ to $N_1$ and $\mu=\mu_{((CTC^{-1})|_{C^{-1}(N_2)})^{*}}=\Gamma_C((P_{N_2}TP_{N_1})^{'}|_{N_1^{'}})$. Therefore, $\Gamma_C$ is fully faithful.  

To show that $\Gamma_C$ is a local isomorphism, we need to prove that for every $N^{'}\in \textit{\textbf v}\mathscr F(H^{'})$, $\Gamma_C|_{\langle N^{'}\rangle}$ is an isomorphism from $\langle N^{'}\rangle$ onto $\langle \Gamma_C(N^{'})\rangle$. Since $\Gamma_C$ is fully faithful, it suffices to see that $\textit{\textbf v}\Gamma_C$ is a bijection from $\textit{\textbf v} \langle N^{'}\rangle$ to $\textit{\textbf v} \langle \Gamma_C(N^{'})\rangle=\textit{\textbf v} \langle H(\gamma^{Q_{C^{-1}(N)}}, -)\rangle$. Suppose $H(\gamma^P, -)\subseteq H(\gamma^{Q_{C^{-1}(N)}}, -)$, then there exists a unique epimorphism $A$ from $C^{-1}(N)$ to $R(P)$ such that $\gamma^P=\gamma^{Q_{C^{-1}(N)}} \star A$. Thus, $P=Q_{C^{-1}(N)}AP_{R(P)}$ implies that $P$ is a bounded operator with $N^\perp\subseteq C(Z(P))$. Therefore, if we let $M=C(Z(P))^\perp$, then $M\subseteq N$ and $\gamma^P=\gamma^{Q_{C^{-1}(M)}} \star P|_{C^{-1}(M)}$. Since, $P|_{C^{-1}(M)}$ is an isomorphism from $C^{-1}(M)$ to $R(P)$, it follows that $H(\gamma^{Q_{C^{-1}(M)}}, -)=H(\gamma^P, -)$, hence, for each $H(\gamma^P, -)\subseteq H(\gamma^{Q_{C^{-1}(N)}}, -)$, there exists a unique subspace $M^{'}$ of $N^{'}$ such that $\Gamma_C(M^{'})=H(\gamma^P, -)$.
In addition, for each $N\in \textit{\textbf v}\mathscr F(H)$, there exists $C(N)^{'}\in \textit{\textbf v}\mathscr F(H^{'})$ such that $N\in \mathcal MH(\gamma^{Q_N}, -)=\mathcal M\Gamma_C(C(N)^{'})$. Therefore, $\Gamma_C:\mathscr F(H^{'})\rightarrow N^*\mathscr F(H)$ is a cross-connection. We summarize the foregoing discussion as follows:
\begin{theorem}
	Let $C$ be an invertible bounded operator on a Hilbert space $H$. Then, $(\mathscr F(H^{'}), \mathscr F(H); \Gamma_C)$ is a cross-connection from $\mathscr F(H^{'})$ to $\mathscr F(H)$.
\end{theorem}

It is seen that the semigroups of normal cones of isomorphic normal categories are also isomorphic \cite{s}. Hence, Theorem 4.2 of \cite{a} shows that $\mathscr F(H)$ and $\mathscr F(H^{'})$ are isomorphic as normal categories, ie., $T\mathscr F(H)\cong T\mathscr F(H^{'})$. So, for each normal cone $\tau$ in $\mathscr F(H^{'})$ with vertex $N^{'}$, there exists a unique operator $T\in \mathcal S^{'}$ with $R(T)=N$ such that for $M^{'} \in \textit{\textbf v}\mathscr F(H^{'}),$
$$\tau(M^{'})=(P_N(P_MT)^*P_M)^{'}|_{M^{'}}: M^{'}\rightarrow N^{'},$$
which maps $f_m\in M^{'}$ to $f_{T(m)}\in N^{'}$, also it is known that for a normal cone $\gamma$ in a normal category $\mathscr C$, $\gamma\in E(T\mathscr C)$ if and only if $\gamma(c_\gamma)=1_{c_\gamma}$ (see \cite[Theorem 2]{s}). Thus, if $\tau$ is an idempotent normal cone in $\mathscr F(H^{'})$ with vertex $N^{'}$ 
and is induced by the finite rank operator $T$ on $H$, then $T|_N$ is the identity operator on $N$ implies that $T\in E(\mathcal S^{'})$.
In addition, there exists an isomorphism from $\mathscr F(H)$ to $\mathscr F(H^{'})$ such that the restriction of this isomorphism from the hom-sets of $\mathscr F(H)$ to the hom-sets of $\mathscr F(H^{'})$ is a surjective conjugate-linear isometry (see \cite[Theorem 4.2]{a}). Hence, the set $B\mathscr F(H^{'})$ of all bounded normal cones in $\mathscr F(H^{'})$ is isomorphic to the set $\mathcal S$ of all finite rank bounded operators on $H$. For $T\in \mathcal S$ with $N=R(T^*)$, let $\tau^T$ denote the bounded normal cone induced by $T^*$ in $\mathscr F(H^{'})$ with vertex $N^{'}$. Specifically, the bounded normal cone $\tau^T$ is defined for $M^{'} \in \textit{\textbf v}\mathscr F(H^{'})$ by
$$\tau^T(M^{'})=(P_N(P_MT^*)^*P_M)^{'}|_{M^{'}}=(P_NTP_M)^{'}|_{M^{'}}: M^{'}\rightarrow N^{'}.$$ 
Moreover, by defining addition, scalar multiplication, norm and multiplication on $B\mathscr F(H^{'})$ by:
\begin{align*}
	\tau^{T_1}+\tau^{T_2}&=\tau^{T_1+T_2},\\
	k\tau^T&=\tau^{kT},\\
	\lVert \tau^T\rVert&=\lVert T\rVert \text{ and}\\
	\tau^{T_1} \cdot \tau^{T_2}&=\tau^{T_2T_1}
\end{align*}
for $T, T_1, T_2 \in \mathcal S$ and $k\in \mathbf K$, by Theorems 4.3 and 5.3 of \cite{a} it is seen that $B\mathscr F(H^{'})$ is a normed algebra and is conjugate-linear isomorphic to the normed algebra $\mathcal S$ and $\tau^T \in E(B\mathscr F(H^{'}))$ if and only if $T \in E(\mathcal S)$.
Let $\tau^P$ be an idempotent bounded normal cone in $\mathscr F(H^{'})$ with vertex $M^{'}$, and  $N=R(P)$. Then, $P=PP_N$ gives $\tau^P=\tau^{P_N} \cdot \tau^P=\tau^{P_N}\star(P_MPP_N)^{'}|_{N^{'}}$, since $P^*|_N$ is an isomorphism from $N$ to $M$, $H(\tau^P, -)=H(\tau^{P_N}, -)$. Thus, we have the following lemma.
\begin{lemma}
	Let $H^{'}$ be the dual space of $H$, and $\mathscr F(H^{'})$ be the category of finite-dimensional subspaces of $H^{'}$. Then, the bounded normal dual $B^*\mathscr F(H^{'})$ of $\mathscr F(H^{'})$ is the full subcategory of the normal dual $N^*\mathscr F(H^{'})$ of $\mathscr F(H^{'})$ such that
$$\textit{\textbf v}B^*\mathscr F(H^{'})=\{H(\tau^{P_N}, -): N\in \textit{\textbf v}\mathscr F(H)\}.$$
\end{lemma}
For $N, M_1, M_2 \in \textit{\textbf v}\mathscr F(H)$, 
\begin{align*}
	H(\tau^{P_N}, M_1^{'})&=\{\tau^{P_N}\star((P_{M_1}AP_N)^{'}|_{N^{'}})^\circ: A\in \mathscr F(H)(M_1, N)\}\\
	&=\{\tau^{P_{M_1}AP_N}:  A\in \mathscr F(H)(M_1, N)\},
\end{align*}
and for $(P_{M_2}TP_{M_1})^{'}|_{M_1^{'}}:M_1^{'}\rightarrow M_2^{'},$
\begin{align*}
	H(\tau^{P_N}, (P_{M_2}TP_{M_1})^{'}|_{M_1^{'}})&:H(\tau^{P_N}, M_1^{'})\rightarrow H(\tau^{P_N}, M_2^{'})\\
	&:\tau^{P_{M_1}AP_N}\mapsto \tau^{P_{M_2}TAP_N}.
\end{align*}
Also, for each natural transformation $\varphi: H(\tau^{P_{N_1}}, -)\rightarrow H(\tau^{P_{N_2}}, -)$, there exists a unique linear map $T:N_1\rightarrow N_2$ such that for $M^{'}\in \textit{\textbf v}\mathscr F(H^{'})$ and $A\in \mathscr F(H)(M, N_1)$,
\begin{align*}
	\varphi(M^{'})&: H(\tau^{P_{N_1}}, M^{'}) \rightarrow H(\tau^{P_{N_2}}, M^{'})\\
	&:\tau^{P_MAP_{N_1}} \mapsto \tau^{P_MATP_{N_2}},
\end{align*}
we denote the above natural transformation determined by $T$ by $\varphi_T$.

For $N, M \in \textit{\textbf v}\mathscr F(H)$, if $A\in \mathscr F(H)(C^{-1}(N), M)$, then $\gamma^{Q_{C^{-1}(N)}AP_M}\in H(\gamma^{Q_{C^{-1}(N)}}, M)=\Gamma_C(N^{'})(M)$ and $((C^{-1}AC)|_N)^*$ is a linear map from $C(M)$ to $N$, hence, $\tau^{P_NC^{-1}ACP_{C(M)}}\in H(\tau^{P_{C(M)}}, N^{'})$. Now, define $\Delta_C:\mathscr F(H) \rightarrow N^{*}\mathscr F(H^{'})$ by 
$$\Delta_C(M)=H(\tau^{P_{C(M)}}, -)$$
for $M \in \textit{\textbf v}\mathscr F(H)$ and for $S: M_1\rightarrow M_2$,
$$\Delta_C(S)=\varphi_{(C^{-1}SC)|_{C(M_1)}}:H(\tau^{P_{C(M_1)}}, -)\rightarrow H(\tau^{P_{C(M_2)}}, -),$$
then, $\Delta_C$ is a local isomorphism from $\mathscr F(H)$ to $N^{*}\mathscr F(H^{'})$. For $S_1: M_1\rightarrow M_2$ and $S_2: M_2\rightarrow M_3$,
\begin{align*}
	\Delta_C(S_1S_2)&=\varphi_{(C^{-1}S_1S_2C)|_{C(M_1)}}\\
	&=\varphi_{(C^{-1}S_1C)|_{C(M_1)}}\varphi_{(C^{-1}S_2C)|_{C(M_2)}}\\
	&=\Delta_C(S_1)\Delta_C(S_2).
\end{align*}
Further, if $M\subseteq N$, then $P_N|_M$ is the inclusion from $M$ to $N$ and $\Delta_C(P_N|_M)=\varphi_{(C^{-1}(P_N|_M)C)|_{C(M)}}$ is the inclusion from $H(\tau^{P_{C(M)}}, -)$ to $H(\tau^{P_{C(N)}}, -)$. Therefore, $\Delta_C$ is an inclusion preserving functor. Also, for each natural transformation $\varphi: H(\tau^{P_{C(M_1)}}, -)\rightarrow H(\tau^{P_{C(M_2)}}, -),$ there exists a unique linear map $T: C(M_1)\rightarrow C(M_2)$ such that $\varphi=\varphi_T.$  If we let $S=(CTC^{-1})|_{M_1}$, then $S$ is a linear from $M_1$ to $M_2$ and $\varphi=\varphi_{(C^{-1}SC)|_{C(M_1)}}=\Delta_C(S).$ That is, $\Delta_C$ is fully faithful.
Finally, to prove that for each $M\in \textit{\textbf v}\mathscr F(H)$, $\textit{\textbf v}\Delta_C$ is a bijection from $\textit{\textbf v}\langle M \rangle$ to $\textit{\textbf v}\langle H(\tau^{P_{C(M)}}, -) \rangle.$ For, consider $\tau \in E(T\mathscr F(H^{'}))$ with vertex $N^{'}$ and $H(\tau, -)\subseteq H(\tau^{P_{C(M)}}, -)$, then there is a unique epimorphism $A^{*}$ from $C(M)$ to $N$ such that $\tau=\tau^{P_{C(M)}}\star(P_NAP_{C(M)})^{'}|_{C(M)^{'}}=\tau^{P_NAP_{C(M)}}$, that is $\tau$ is induced by $(P_NAP_{C(M)})^{*}\in \mathcal S$. However, since $\tau\in E(B\mathscr F(H^{'}))$, there exists a unique idempotent $P\in \mathcal S$ such that $P=P_NAP_{C(M)}$. Let $M_1=R(P)$, then $M_1\subseteq C(M)$ and $\tau^P=\tau^{P_{M_1}} \star (P_NPP_{M_1})^{'}|_{M_1^{'}}$. We have $P^{*}|_{M_1}$ is an isomorphism from $M_1$ to $N$, which implies that $H(\tau^{P_{M_1}}, -)=H(\tau, -)$. Hence, for each $H(\tau, -)\subseteq H(\tau^{P_{C(M)}}, -)$, there exists a unique subspace $C^{-1}(M_1)$ of $M$ such that $\Delta_C(C^{-1}(M_1))=H(\tau, -)$. Thus, $\Delta_C$ is a local isomorphism.

Moreover, for each $N^{'}\in \textit{\textbf v}\mathscr F(H^{'})$, there exists $C^{-1}(N)\in \textit{\textbf v}\mathscr F(H)$ such that $N^{'}\in \mathcal M H(\tau^{P_N}, -)=\mathcal M\Delta_C(C^{-1}(N))$. Therefore, $(\mathscr F(H), \mathscr F(H^{'}); \Delta_C)$ is a cross-connection, as stated in the following theorem:

\begin{theorem}
	Let $C$ be an invertible bounded operator on $H$. Then, the functor $\Delta_C:\mathscr F(H)\rightarrow N^*\mathscr F(H^{'})$ is a cross-connection from $\mathscr F(H)$ to $\mathscr F(H^{'})$.
\end{theorem}
Now, consider the bifunctors $\Gamma_C(-, -):\mathscr F(H)\times \mathscr F(H^{'})\rightarrow \mathbf{Set}$ and $\Delta_C(-, -):\mathscr F(H)\times \mathscr F(H^{'})\rightarrow \mathbf{Set}$ corresponding to the cross-connections $\Gamma_C:\mathscr F(H^{'})\rightarrow N^*\mathscr F(H)$ and $\Delta_C: \mathscr F(H) \rightarrow N^*\mathscr F(H^{'})$. For $(M_1, N_1^{'}), (M_2, N_2^{'})\in \textit{\textbf v}\mathscr F(H)\times \textit{\textbf v}\mathscr F(H^{'})$ and  $(S, (P_{N_2}TP_{N_1})^{'}|_{N_1^{'}}): (M_1, N_1^{'})\rightarrow (M_2, N_2^{'})$,
\begin{align*}
	\Gamma_C(M_1, N_1^{'})&=H(\gamma^{Q_{C^{-1}(N_1)}}, M_1)\\
	&=\{\gamma^{Q_{C^{-1}(N_1)}AP_{M_1}}: A\in \mathscr F(H)(C^{-1}(N_1), M_1)\},\\
	\Gamma_C(S, (P_{N_2}TP_{N_1})^{'}|_{N_1^{'}})&:H(\gamma^{Q_{C^{-1}(N_1)}}, M_1)\rightarrow H(\gamma^{Q_{C^{-1}(N_2)}}, M_2)\\
	&:\gamma^{Q_{C^{-1}(N_1)}AP_{M_1}}\mapsto \gamma^{Q_{C^{-1}(N_2)}CTC^{-1}ASP_{M_2}},\\
	\Delta_C(M_1, N_1^{'})&=H(\tau^{P_{C(M_1)}}, N_1^{'})\\
	&=\{\tau^{P_{N_1}BP_{C(M_1)}}: B\in \mathscr F(H)(N_1, C(M_1))\}\text{ and}\\
	\Delta_C(S, (P_{N_2}TP_{N_1})^{'}|_{N_1^{'}})&:H(\tau^{P_{C(M_1)}}, N_1^{'})\rightarrow H(\tau^{P_{C(M_2)}}, N_2^{'})\\
	&:\tau^{P_{N_1}BP_{C(M_1)}}\mapsto \tau^{P_{N_2}TBC^{-1}SCP_{C(M_2)}}.
\end{align*}
Now define $\chi_C:\Gamma_C(-, -) \rightarrow \Delta_C(-, -)$ for $(M_1, N_1^{'}) \in \textit{\textbf v}\mathscr F(H)\times \textit{\textbf v}\mathscr F(H^{'})$ by
\begin{align*}
\chi_C(M_1, N_1^{'})&:\Gamma_C(M_1, N_1^{'}) \rightarrow \Delta_C(M_1, N_1^{'})\\
&:\gamma^{Q_{C^{-1}(N_1)}AP_{M_1}}\mapsto \tau^{P_{N_1}C^{-1}ACP_{C(M_1)}},
\end{align*}
for $(M_1, N_1^{'}), (M_2, N_2^{'})\in \textit{\textbf v}\mathscr F(H)\times \textit{\textbf v}\mathscr F(H^{'})$ and $(S, (P_{N_2}TP_{N_1})^{'}|_{N_1^{'}}): (M_1, N_1^{'})\rightarrow (M_2, N_2^{'})$ the following diagram commutes:

\[
\xymatrix@=2.5cm{\Gamma_C(M_1, N_1^{'})\ar[d]_{\Gamma_C\big(S, (P_{N_2}TP_{N_1})^{'}|_{N_1^{'}}\big)}\ar[r]^{\chi_C(M_1, N_1^{'})} & \Delta_C(M_1, N_1^{'})\ar[d]^{\Delta_C\big(S, (P_{N_2}TP_{N_1})^{'}|_{N_1^{'}}\big)}\\ 
	 \Gamma_C(M_2, N_2^{'})\ar[r]_{\chi_C (M_2, N_2^{'})}  & \Delta_C (M_2, N_2^{'})
}
\]
Therefore, $\chi_C$ is a natural transformation from $\Gamma_C(-, -)$ to $\Delta_C(-, -)$, and  for each normal cone $\tau=\tau^{P_{N_1}BP_{C(M_1)}}\in \Delta_C(M_1, N_1^{'})$, there exists a unique normal cone $\gamma=\gamma^{Q_{C^{-1}(N_1)}CBC^{-1}P_{M_1}}\in \Gamma_C(M_1, N_1^{'})$ such that $\chi_C(M_1, N_1^{'})(\gamma)=\tau.$ That is, for each $(M, N^{'}) \in \textit{\textbf v}\mathscr F(H)\times \textit{\textbf v}\mathscr F(H^{'})$, $\chi_C(M, N^{'})$ is a bijection from $\Gamma_C(M, N^{'})$ to $\Delta_C(M, N^{'})$. Hence, $\chi_C$ is a natural isomorphism from $\Gamma_C(-, -)$ to $\Delta_C(-, -)$. Consequently, we obtain the following theorem.
\begin{theorem}
Let $C$ be an invertible bounded operator on $H$. Then, $\Delta_C$ is the dual cross-connection of $\Gamma_C$, and  $\chi_C$ is the duality associated with $\Gamma_C$.
\end{theorem}
For $T\in \mathcal S$, if $N\supseteq C(Z(T))^\perp$ and $M\supseteq R(T)$, then $\gamma^T=\gamma^{Q_{C^{-1}(N)}TP_{M}}\in \Gamma_C(M, N^{'})$ and $\chi_C(M, N^{'})(\gamma^{T})=\tau^{C^{-1}TC}=\tau^{P_NC^{-1}TCP_{C(M)}}\in \Delta_C(M, N^{'}).$ Therefore, $T \mapsto (\gamma^{T}, \tau^{C^{-1}TC})$ is a one-to-one correspondence from the set $\mathcal S$ to the set $\tilde{S}\Gamma_C$ of pairs of linked cones. Denote this correspondence by $\eta$, i.e., $\eta:\mathcal S \rightarrow \tilde{S}\Gamma_C$ is defined for $T\in\mathcal S$ by 
$$\eta(T)=(\gamma^T, \tau^{C^{-1}TC}).$$
Let $T_1\in\mathcal S$ with $N_1 \supseteq C(Z(T_1))^\perp$ and $M_1 \supseteq R(T_1)$, and $T_2\in\mathcal S$ with $N_2 \supseteq C(Z(T_2))^\perp$ and $M_2 \supseteq R(T_2)$. Then, $\gamma^{T_1T_2} \in \Gamma_C(M_2, N_1^{'})$ and
$$(\gamma^{T_1T_2}, \tau^{C^{-1}T_1T_2C})=(\gamma^{T_1}\cdot \gamma^{T_2}, \tau^{C^{-1}T_2C} \cdot \tau^{C^{-1}T_1C})$$ is a pair of linked cones. Also, the multiplication on $\tilde{S}\Gamma_C$ defined by
$$(\gamma_1, \tau_1)(\gamma_2, \tau_2)=(\gamma_1\cdot\gamma_2, \tau_2\cdot \tau_1)$$
for $(\gamma_1, \tau_1), (\gamma_2, \tau_2) \in \tilde{S}\Gamma_C$, satisfies  $\eta(T_1T_2)=\eta(T_1)\eta(T_2)$, hence, $\tilde{S}\Gamma_C$ is a semigroup isomorphic to the regular semigroup $\mathcal S$ and for $N \supseteq C(Z(T_1+T_2))^\perp,\,M \supseteq R(T_1+T_2)$, $\gamma^{T_1+T_2}\in \Gamma_C(M, N^{'})$, 
$$(\gamma^{T_1+T_2}, \tau^{C^{-1}(T_1+T_2)C})=(\gamma^{T_1}+ \gamma^{T_2}, \tau^{C^{-1}T_1C}+ \tau^{C^{-1}T_2C})$$ is a pair of linked cones. Thus, for $(\gamma_1, \tau_1), (\gamma_2, \tau_2) \in \tilde{S}\Gamma_C$, 
$$(\gamma_1, \tau_1)+(\gamma_2, \tau_2)=(\gamma_1+\gamma_2, \tau_1+\tau_2)$$
defines an addition on $\tilde{S}\Gamma_C$, and it follows that $\eta(T_1+T_2)=\eta(T_1)+\eta(T_2)$. Also, for $k\in \mathbf K$, $\gamma^{kT_1} \in \Gamma_C(M_1, N_1^{'})$ and 
$$(\gamma^{kT_1}, \tau^{C^{-1}(kT_1)C})=(k\gamma^{T_1}, k\tau^{C^{-1}T_1C})$$
is a pair of linked cones. i.e.,
$$k(\gamma, \tau)=(k\gamma, k\tau)$$
for $(\gamma, \tau)\in \tilde{S}\Gamma_C$ defines a scalar multiplication on $\tilde{S}\Gamma_C$ and $\eta(kT_1)=k\eta(T_1)$. Thus, the following proposition:
\begin{proposition}
	Let $C$ be an invertible bounded operator on $H$. Then, the pairs of linked cones
	$\tilde{S}\Gamma_C$ is a regular semigroup, moreover it admits a vector space structure and is isomorphic to the vector space $\mathcal S$ of all finite rank bounded operators on $H$.
\end{proposition}
Also it is seen that for $(\gamma, \tau)\in \tilde{S}\Gamma_C$, 
$$\lVert (\gamma, \tau) \rVert =\max\{\lVert \gamma \rVert, \lVert \tau \rVert\}$$
is a norm on $\tilde{S}\Gamma_C$ and so $\tilde{S}\Gamma_C$ is a normed algebra. If $C$ is a unitary operator on $H$, then for $T\in\mathcal S$,
\begin{align*}
	\lVert \gamma^{T} \rVert&= \lVert T \rVert\text{ and}\\
	\lVert \tau^{C^{*}TC} \rVert &=\lVert C^{*}TC \rVert =\lVert T \rVert.
\end{align*}
 That is, for $T\in \mathcal S$, $\lVert \eta(T) \rVert=\lVert (\gamma^{T}, \tau^{C^*TC}) \rVert=\lVert T \rVert$ when $C$ is a unitary operator on $H$. Hence, we obtain the following theorem.
\begin{theorem}
	Let $C$ be a unitary operator on a Hilbert space $H$. Then, the regular semigroup $\tilde{S}\Gamma_C$ of pairs of linked cones is a normed algebra isomorphic to the normed algebra $\mathcal S$ of all finite rank bounded operators on $H$.
\end{theorem}	

\section{Compact operators}\label{sec4}

Let $\mathscr P(H)$ denote the category of proper subspaces of $H$, whose vertices are proper subspaces of $H$, and for $M, N\in \textit{\textbf v} \mathscr P(H)$, the morphisms from $M$ to $N$ are linear maps from $M$ to $N$. Then, the category $\mathscr P(H)$ with usual subspace inclusions as inclusions is a category with subobjects. Similar to $\mathscr F(H)$, it is seen that $\mathscr P(H)$ is a normal category \cite{a}. By modifying the proof of Theorem 5.1 in \cite{a}, it follows that the semigroup $T\mathscr P(H)$ of all normal cones in $\mathscr P(H)$ is isomorphic to the semigroup of non-surjective operators on $H$. Precisely, for each normal cone $\delta$ in $\mathscr P(H)$ with vertex $N$, there is a unique non-surjective operator $T$ on $H$ such that $R(T)=N$ and $\delta=\delta^T:\textit{\textbf v} \mathscr P(H)\rightarrow \mathscr P(H)$ is given for $M\in \textit{\textbf v} \mathscr P(H)$ by
$$\delta^T(M)=T|_M:M\rightarrow N.$$

For a finite-dimensional Hilbert space $H$, the spectrum $\sigma(T)$ of an operator $T$ on $H$ is given by
\begin{align*}
    \sigma(T)&=\{\lambda\in \mathbf K: T-\lambda I \text{ is not surjective}\}\\
    &=\{\lambda\in \mathbf K: \delta^{T-\lambda I} \in T\mathscr P(H)\}.
\end{align*}
However, when $H$ be an infinite-dimensional Hilbert space and $\mathcal K$, the Banach algebra of all compact operators on $H$, we have that if $T\in \mathcal K$, then the range space $R(T^*-\bar \lambda I)$ of $T^*-\bar \lambda I$ is a closed subspace of $H$ for $\lambda \neq 0$. Thus, $\lambda \neq 0$ is an eigenvalue of $T\in \mathcal K$ if and only if 
$T^*-\bar \lambda I$ is not surjective. Hence, we have the following:
\begin{theorem}
    Let $H$ be an infinite-dimensional Hilbert space and $T$ be a compact operator on $H$. Then, the set of non-zero eigenvalues of $T$ is
    $$\sigma_p(T)\backslash \{0\}=\{\lambda\in \mathbf K\backslash \{0\}:\delta^{T^*-\bar \lambda I}\in T\mathscr P(H)\},$$ where $\sigma_p(T)$ denote the point spectrum of $T$.
\end{theorem}
Also, since compact operators on an infinite-dimensional Hilbert space are not surjective, $\sigma(T)=\{0\}\cup \sigma_p(T)$ for $T\in \mathcal K$ leads to the following theorem.
\begin{theorem}
Let $H$ be an infinite-dimensional Hilbert space and $T$ be a compact operator on $H$. Then, the spectrum $\sigma(T)$ of $T$ is given by $$\sigma(T)=\{\lambda\in \mathbf K: \delta^{T^*-\bar \lambda I}\in T\mathscr P(H)\}.$$ 
\end{theorem}
Further, for the self-adjoint compact operator $T$ on $H$,  the eigenspace $E_\lambda$ of $T$ corresponding to the eigenvalue $\lambda$ is the orthogonal complement of the vertex of $\delta^{T-\lambda I}$, that is, $E_\lambda=(c_{\delta^{T-\lambda I}})^\perp$. Consequently, the spectral theorem for compact self-adjoint operators on $H$ can be stated in terms of normal cones in $\mathscr P(H)$ as follows:
\begin{theorem}
    Let $T$ be a self-adjoint compact operator on a Hilbert space $H$, and for $\delta^{T-\lambda I}\in T\mathscr P(H)$, let $\mathscr B_\lambda$ be an orthonormal basis for $(c_{\delta^{T-\lambda I}})^\perp$. Then, $\bigcup\limits_{\lambda\in\sigma(T)\backslash\{0\}}\mathscr B_\lambda$ is an orthonormal basis for $Z(T)^\perp$, and $\bigcup\limits_{\lambda\in\sigma(T)}\mathscr B_\lambda$ is an orthonormal basis for $H$, consisting of eigenvectors of $T$.
\end{theorem}
Next, consider an infinite-dimensional Hilbert space $H$, in the light of Theorem 5.2 in \cite{a}, it is seen that the set $B\mathscr P(H)$ of all bounded normal cones in $\mathscr P(H)$ is a semigroup isomorphic to the semigroup of all non-surjective bounded operators on $H$, and each compact operator $T$ on $H$ induces a normal cone $\delta^T$ in $\mathscr P(H)$ with vertex $R(T)$. Since, the Banach algebra $\mathcal K$ of all compact operators on $H$ is not a regular semigroup with respect to function composition, these normal cones $\delta^T$ fails to yield a regular semigroup. However, it is known that for each compact operator $T$ on $H$, there exists a sequence $(T_n)$ in $\mathcal S$ such that $\lVert T_n-T\rVert\rightarrow 0$. Therefore, we define the convergence of sequences of bounded normal cones in $\mathscr P(H)$ as below:
\begin{definition}
	Let $B\mathscr P(H)$ be the semigroup of all bounded normal cones in $\mathscr P(H)$. Then, a sequence $(\delta^{T_n})$ in $B\mathscr P(H)$ is said to converge in $B\mathscr P(H)$ if there exists some $\delta^T\in B\mathscr P(H)$ such that $\lVert T_n-T\rVert\rightarrow 0$. 
\end{definition}
Thus, for each $T\in \mathcal K$, there exists a sequence $(\delta^{T_n})$ in $B\mathscr P(H)$ such that $T_n\in \mathcal S$ and $\delta^{T_n}$ converges to $\delta^T$.  Note that for $T\in \mathcal S$, the restriction of the bounded normal cone $\delta^{T}$ to $\textit{\textbf v}\mathscr F(H)$ is $\gamma^{T}$. Hence, we define addition, scalar multiplication and norm on the subsemigroup $F\mathscr P(H)=\{\delta^T\in B\mathscr P(H): T\in \mathcal S\}$ of $B\mathscr P(H)$ by
\begin{align*}
	\delta^{T_1}+\delta^{T_2}&=\delta^{T_1+T_2},\\
	k\delta^T&=\delta^{kT} \text{ and}\\
	\lVert \delta^T\rVert&=\lVert T\rVert 
\end{align*} 
for $T, T_1, T_2 \in \mathcal S$ and $k\in \mathbf K$. Then, using the map $\phi: B\mathscr F(H) \rightarrow F\mathscr P(H)$ defined for $\gamma^T \in B\mathscr F(H)$ by $\phi(\gamma^T)=\delta^T$, it can be observed that $F\mathscr P(H)$ is a normed algebra and is isomorphic to the normed algebra $B\mathscr F(H)$. Since the closure of a normed algebra is a Banach algebra and the closure $\overline{\mathcal S}$ of $\mathcal S$ is $\mathcal K$, the following theorem is established:
\begin{theorem}
	Let $F\mathscr P(H)$ be the normed algebra contained in $B\mathscr P(H)$ such that $F\mathscr P(H)$ is isomorphic to the normed algebra $B\mathscr F(H)$. Then, the closure $\overline{F\mathscr P(H)}$ of $F\mathscr P(H)$ in $B\mathscr P(H)$ is a Banach algebra isomorphic to the Banach algebra $\mathcal K$ of all compact operators on $H$. 
\end{theorem}



\end{document}